\documentclass[11pt]{amsart}
\title{Selective screenability and the Hurewicz property
\footnote{\lowercase{{\bf {\uppercase{K}}ey words and phrases:} {\uppercase{H}}aver property, selective screenability, totally bounded, $\sigma$-totally bounded, {\uppercase{H}}urewicz property, {\uppercase{M}}enger property, selection principle.\\
{{\bf{\uppercase{S}}ubject {\uppercase{C}}lassification:} {\uppercase{P}}rimary 54{\uppercase{D}}20, 54{\uppercase{D}}45, 55{\uppercase{M}}10; {\uppercase{S}}econdary 03{\uppercase{E}}20.}
}}
}
\author{Liljana Babinkostova}
\date{}
\usepackage{amssymb}
\usepackage{amsfonts}
\newcommand{\op}{\mathcal{O}}

\newcommand{\Sc}{{\sf S}_c}

\newcommand{\naturals}{{\mathbb N}}

\newtheorem{theorem}{Theorem}
\newtheorem{lemma}[theorem]{Lemma}	
\newtheorem{corollary}[theorem]{Corollary}	
\newtheorem{problem}{Problem}
\newtheorem{conjecture}[problem]{Conjecture}
\begin{document}
\begin{abstract} 
We characterize the Hurewicz covering property in metrizable spaces in terms of properties of the metrics of the space - Theorem \ref{hurewicztotbdd}. Then we show that a weak version of selective screenability, when combined with the Hurewicz property, implies selective screenability - Theorem \ref{hurewiczscfin}.
\end{abstract}
\maketitle

\section{Definitions and notation}

Let $X$ be an infinite set, and let $\mathcal{A}$ and $\mathcal{B}$ be collections of families of subsets of $X$.
The selection principle $\Sc(\mathcal{A},\mathcal{B})$, introduced in \cite{lbsc}, states:
\begin{quote}
For each sequence $(A_n:n<\infty)$ of elements of the family $\mathcal{A}$ there exists a sequence $(B_n:n<\infty)$ such that for each $n$ $B_n$ is a pairwise disjoint family refining $A_n$, and $\bigcup_{n<\infty}B_n$ is a member of the family $\mathcal{B}$.
\end{quote}
For $X$ topological space $\op$ denotes the collection of all open covers of $X$ and $\op_{fin}$ denotes the collection of all finite open covers of $X$. For a positive integer $n$ let $\op_n$ denote the collection of open covers consisting of at most $n$ sets. Addis and Gresham introduced the instance $\Sc(\op,\op)$ of the selection principle in \cite{ag}, where it was called property C. It is a selective version of the screenability property introduced by Bing in \cite{bing}. 

As was shown in \cite{ag}, $\Sc(\op,\op)$ is a natural generalization of finite covering dimension to the infinite. Alexandroff's notion of weakly infinite dimensional is also a natural generalization of finite covering dimension, and is equivalent to $\Sc(\op_2,\op)$. Hurewicz's notion of countable dimensionality is another natural generalization of finite covering dimension: $X$ is countable dimensional if it is a union of countably may finite dimensional subspaces. The following implications hold - see \cite{ag}:
\[
   \mbox{countable dimensional} \Rightarrow \Sc(\op,\op)\Rightarrow \Sc(\op_{fin},\op) \Rightarrow \Sc(\op_2,\op).
\]
The Hilbert cube, [0,1]$^{\naturals}$, does not have property $\Sc(\op_2,\op)$ - \cite{ag}. Borst proved in \cite{pborst2} that there exists a compact separable metric space $X$ which has property $\Sc(\op_2,\op)$, but not property $\Sc(\op,\op)$. Since for compact spaces $\Sc(\op_{fin},\op) \Leftrightarrow \Sc(\op,\op)$, Borst's example shows that $\Sc(\op_2,\op)$ does not imply $\Sc(\op_{fin},\op)$. R. Pol gave in \cite{rpol} a compact metric space which has property $\Sc(\op,\op)$ but is not countable dimensional. It is an open problem if $\Sc(\op_{fin},\op)$ implies $\Sc(\op,\op)$ - see Question 3.10 of \cite{borst}. We expect that the answer to this question is ``No", and state a conjecture about it near the end of this paper. In \cite{borst} a class of spaces which do not distinguish $\Sc(\op_{fin},\op)$ and $\Sc(\op,\op)$ is identified. In this paper we will extend this to a larger class of separable metric spaces which do not distinguish $\Sc(\op_{fin},\op)$ and $\Sc(\op,\op)$. Examples show that the class we describe properly extends the class from \cite{borst}.

In Section 2 we first give a convenient characterization of the Hurewicz property in metrizable spaces. In Section 3 we show that metrizable spaces with the Hurewicz property do not distinguish $\Sc(\op_{fin},\op)$ and $\Sc(\op,\op)$. In Section 4 we connect this with Borst's work from \cite{borst} and in the final section we state a conjecture. 

\section{Characterizing the Hurewicz property in metrizable spaces.}

A topological space $X$ has the \emph{Hurewicz property} \cite{wh} if there is for each sequence $(\mathcal{U}_n:n<\infty)$ of open covers of $X$ a sequence $(\mathcal{V}_n:n<\infty)$ such that for each $n$, $\mathcal{V}_n$ is a finite subset of $\mathcal{U}_n$, and each element of $X$ is in all but finitely many of the sets $\cup\mathcal{V}_n$. The metrizable space $X$ is said to be \emph{Haver} \cite{h} with respect to a metric $d$ if there is for each sequence $(\epsilon_n:n<\infty)$ of positive reals a sequence $(\mathcal{V}_n:n<\infty)$ where each $\mathcal{V}_n$ is a pairwise disjoint family of open sets, each of $d$-diameter less than $\epsilon_n$, such that $\bigcup_{n<\infty}\mathcal{V}_n$ is a cover of $X$. 

A metric space $(X,d)$ is totally bounded if there is for each $\epsilon>0$ a finite set $F\subset X$ such that $X\subseteq\bigcup_{f\in F}B_d(f,\epsilon)$, where $B_d(f,\epsilon)=\{x\in X:\,d(x,f)<\epsilon\}$.
A metric space is $\sigma$-totally bounded if it is a union of countably many subsets, each totally bounded.\\

\begin{theorem}\label{hurewicztotbdd}
Let $(X,d)$ be a metrizable space. The following are equivalent:
\begin{enumerate}
\item{$X$ has the Hurewicz property.}
\item{$X$ is $\sigma$-totally bounded in each equivalent metric.}
\end{enumerate}
\end{theorem}
{\bf Proof:} 
$1\Rightarrow 2$: For each $n$ let ${\delta}_n=(1/2)^{2^n}$ and $\mathcal{U}_n=\{B_d(x,\delta_n):x\in X\}$ where $d$ is an arbitrary fixed metric of $X$. Apply the Hurewicz property to $(\mathcal{U}_n:n<\infty)$. For each $n$ choose a finite set $\mathcal{V}_n\subset\mathcal{U}_n$ such that each $x\in X$ is in all but finitely many of the sets $\cup\mathcal{V}_n$. For each $n$ define $X_n=\bigcap_{m\geq n}\cup\mathcal{V}_m $. Then for each $n$, and for $m\le n$, $X_m\subseteq X_n$ and $\bigcup_{n<\infty} X_n$ covers $X$. We show that each $X_n$ is totally bounded in the metric $d$: Consider an $\epsilon>0$, and consider any $X_n$. Choose $m>n$ so large that $(1/2)^{2^m}\le\epsilon$. Each element of $\mathcal{V}_m$ is an open set of diameter less than $(1/2)^{2^m}$, and $\mathcal{V}_m$  is a finite cover of $X_n$.\\
$2\Rightarrow 1$: Let  $(\mathcal{U}_n:n<\infty)$ be a sequence of open covers of $X$. By Remark 4, page 196 from \cite{JD} let $d$ be a metric generating the topology of $X$ such that for each $n$, $\mathcal{W}_n=\{B_d(x,1/n):x\in X\}$ refines  $\mathcal{U}_n$. Write $X=\bigcup_{n<\infty} X_n$, where each $X_n$ is totally bounded. Choose for each $m$ a finite $\mathcal{F}_m\subset\mathcal{W}_m$ with $X_m\subseteq\cup\mathcal{F}_m$. Then, for each $m$ choose a finite  $\mathcal{V}_m\subset\mathcal{U}_m$ such that $\mathcal{F}_m$ refines $\mathcal{V}_m$. Then, for each $x\in X$ for all but finitely many $n$, $x\in\cup\mathcal{V}_n$.
$\diamondsuit$

\section{$\Sc(\op_{fin},\op)$ in metrizable spaces with the Hurewicz property.}

For easy reference, we denote the following strong form of $\Sc(\op_{fin},\op)$ by the symbol $\Sc^+(\op_{fin},\op)$:
\begin{quote}
For each sequence $(\mathcal{U}_n:n<\infty)$ of finite open covers of $X$ there are a sequence $(\mathcal{W}_n:n<\infty)$ and a sequence $m_1<m_2<...<m_k<...$ such that 
\begin{enumerate}
\item{each $\mathcal{W}_n$ is a finite pairwise disjoint family of open sets,}
\item{each $\mathcal{W}_n$ refines $\mathcal{U}_n$ and}
\item{ for each $x\in X$, for all but finitely many $k$ there is a $j\in[m_k,m_{k+1})$ with $x\in\cup\mathcal{W}_j$.}
\end{enumerate}
\end{quote} 

\begin{lemma}
Let $(X,d)$ be a metrizable space. If $X$ has $\Sc(\mathcal{O}_{fin},\mathcal{O})$ and the Hurewicz property then it has the property $\Sc^+(\op_{fin},\op)$. 
\end{lemma}
{\bf Proof:} 
Recall that $X$ has the Hurewicz property if and only if ONE has no winning strategy in the Hurewicz game (Theorem 27 of \cite{coc1}). Let $(\mathcal{U}_n:n<\infty)$ be a sequence of finite open covers of $X$. Applying $\Sc(\mathcal{O}_{fin},\mathcal{O})$ to $(\mathcal{U}_n:n<\infty)$, choose for each $n$ a pairwise disjoint refinement $\mathcal{V}_n$ of $\mathcal{U}_n$ so that $F(\emptyset)=\bigcup_{n<\infty}\mathcal{V}_n$ covers $X$. This defines ONE's first move in the Hurewicz game. When TWO chooses a finite $T_1\subset F(\emptyset)$, define $m_1=min\{n: T_1\subseteq\bigcup_{j< n}\mathcal{V}_n\}$. Next, apply $\Sc(\mathcal{O}_{fin},\mathcal{O})$ to $(\mathcal{U}_n:n\ge m_1)$ and choose for each $n\ge m_1$ a pairwise disjoint $\mathcal{V}_n$ that refines $\mathcal{U}_n$ consisting of open sets, so that $F(T_1)=\bigcup_{n>m_1}\mathcal{V}_n$ covers $X$. This defines ONE's response to TWO's move $T_1$. When TWO chooses $T_2\subset F(T_1)$, define $m_2=min\{n: T_2\subseteq\bigcup_{m_1\le j <n}\mathcal{V}_j\}$ and apply $\Sc(\mathcal{O}_{fin},\mathcal{O})$ to $(\mathcal{U}_n:n\ge m_2)$ to define $F(T_1,T_2)$, and so on. 

Since $X$ has the Hurewicz property $F$ is not a winning strategy for ONE. Consider an $F$-play $F(\emptyset),T_1,F(T_1),T_2,F(T_1,T_2),T_3...$ lost by ONE. Then each $T_m$ is finite and each $x\in\cup T_m$  for all but finitely many $m$. For $j< m_1$ define $\mathcal{W}_j=\{T\in T_1: (\exists U\in \mathcal{U}_j)(T\subseteq U)\}$. For $m_k\le j < m_{k+1}$ define $\mathcal{W}_j=\{T\in T_{k+1}: (\exists U\in \mathcal{U}_j)(T\subseteq U)\}$. For each $j$, $\mathcal{W}_j$ is finite pairwise disjoint and refines $\mathcal{U}_j$. $\diamondsuit$

\begin{theorem}\label{scfinandhaver}
If $(X,d)$ is $\sigma$-totally bounded and has property $\Sc^+(\op_{fin},\op)$, then $X$ has the Haver property in $d$.
\end{theorem}
{\bf Proof:} 
Write $X=\bigcup_{n<\infty} X_n$, where each $X_n\subset X$ is $d$-totally bounded and $X_n\subset X_{n+1}$. Let $(\epsilon_n: n<\infty)$ be a sequence of positive reals. By replacing $\epsilon_n$'s if necessary, we may assume that always $\epsilon_{n+1}<\frac{1}{2}\cdot\epsilon_n$. For each $n$, put $\delta_n=\frac{2^{2^n}-1}{2^{2^n}}\cdot(\frac{1}{2}\cdot\epsilon _n)$. For each $n$, choose a finite set $F_n\subset X_n$ such that $\{B(x,\delta_n): x\in F_n\}$ covers $X_n$, and put $\mathcal{U}_n = \{B(x,\frac{1}{2}\cdot\epsilon_n): x\in F_n\}\bigcup\{X\setminus\bigcup\{\overline{B(x,\delta_n)}:x\in F_n\}$, a finite open cover of $X$. Observe that for each $n$, $\overline{B(x,\delta_n)}\subset B(x,\epsilon_n)$, and $X_n\bigcap (X\setminus\bigcup\{\overline{B(x,\delta_n)}:x\in F_n\})=\emptyset$.

Apply $\Sc^+(\op_{fin},\op)$ to the sequence $(\mathcal{U}_n:n<\infty)$. For each $n$ find a finite pairwise disjoint refinement $\mathcal{H}'_n$ of $\mathcal{U}_n$ and find a sequence $m_1<m_2<...<m_k<...$ such that for each $x\in X$ for all but finitely many $k$, there is a $j$ with $m_k\le j< m_{k+1}$ and $x\in \cup\mathcal{H}'_j$. Now for each $n$, put 
\[
  \mathcal{H}_n=\{V\in \mathcal{H}'_n: (\exists x\in F_n)(V\subseteq B(x,\frac{1}{2}\cdot\epsilon_n))\}.
\] 
{\bf Claim}: $\bigcup_{n<\infty}\mathcal{H}_n$ covers $X$.\\
For consider $x\in X$. Choose $N$ so large so that for all $n\geq N$, $x\in X_n$ and for all $m_k\ge N$, there is $j\in [m_k,m_{k+1})$ with $x\in\cup\mathcal{H}'_j$. Choose $k$ with $m_k\geq N$ and $j$ with $m_k\leq j<m_{k+1}$ with $x\in V$ for some $V\in\mathcal{H}'_j$. We have that $x\in X_j$, so $V$ is not a subset of $X\setminus (\bigcup\{\overline{B(y,\delta_j)}:y\in F_j\})$ which means that $V\in\mathcal{H}_j$.  

Since the diameter of any element of an $\mathcal{H}_n$ is less than $\epsilon_n$, the sequence $(\mathcal{H}_n:n<\infty)$ witnesses the Haver property of $X$ for $(\epsilon_n: n<\infty)$.
$\diamondsuit$

Note that the Hurewicz property plus $\Sc(\op_2,\op)$ does not imply the Haver property: For if this were to imply the Haver property, then by Theorem 1 of \cite{lbhaver} it would follow that $\Sc(\op_2,\op)$ plus the Hurewicz property implies $\Sc(\op,\op)$. Compactness implies the Hurewicz property, and \cite{pborst2} shows that $\Sc(\op_2,\op)$ plus compact does not imply $\Sc(\op,\op)$. 

\begin{theorem}\label{hurewiczscfin}
If $X$ is a metrizable space and has the Hurewicz property, then the following are equivalent:
\begin{enumerate}
\item{$X$ has $\Sc(\mathcal{O},\mathcal{O})$}
\item{$X$ has $\Sc(\mathcal{O}_{fin},\mathcal{O})$}
\end{enumerate}
\end{theorem}

{\bf Proof:} 
$1\Rightarrow 2$: It is clear.\\
$2\Rightarrow 1$: By the previous theorem $X$ has the Haver property. By Theorem 1 from \cite{lbhaver}we have that $X$ has $\Sc(\mathcal{O},\mathcal{O})$.
$\diamondsuit$

\section{An extension of the class of ``finite C-spaces".}

In \S 3 of \cite{borst}, Borst introduces the notion of a ``finite C-space": A topological space $X$ is a \emph{finite C-space} if there is for each sequence $(\mathcal{U}_n:n<\infty)$ of finite open covers of $X$ an $n$, and a sequence $(\mathcal{V}_j:j\le n)$ such that each $\mathcal{V}_j$ is a disjoint refinement of $\mathcal{U}_j$, and $\bigcup_{j\le n}\mathcal{V}_j$ is an open cover of $X$. And a space $X$ is said to have ``property K" if it has a compact subset $C$ such that for every open subset $U$ of $X$ with $C\subset U$, the set $X\setminus U$ is finite dimensional. And in Theorem 3.8 of \cite{borst} the following equivalence is proved:
\begin{theorem}[Borst]\label{borstfinC} For separable metric spaces $X$ the following are equivalent:
\begin{enumerate}
\item{$X$ is a finite C-space.}
\item{$X$ has $\Sc(\op,\op)$ and property K.}
\end{enumerate}
\end{theorem}
Thus, also in the class of spaces with property K, $\Sc(\op,\op)$ is equivalent to $\Sc(\op_{fin},\op)$. And there are spaces with property K and $\Sc(\op,\op)$ which do not have the Hurewicz property: Let $C$ be the the compact metric space from \cite{rpol}: It has property $\Sc(\op,\op)$ and is infinite dimensional. Let $P$ be the space of irrational numbers. Then $X$, the topological sum of $C$ and $P$, has $\Sc(\op,\op)$ and property K. It is well known that the closed subset P of $X$ does not have the Hurewicz property, and so $X$ does not have the Hurewicz property.

As pointed out in \cite{borst}, the space ${\sf K}_{\omega}$ consisting of the elements $x$ of [0,1]$^\naturals$ for which $x(n)>0$ for only finitely many $n$ is not a ``finite C-space": For if it were a finite C-space, then by Theorem 1.2 of \cite{borst} it has a compactification with property $\Sc(\op,\op)$.  But no compactification of ${\sf K}_{\omega}$ has the property $\Sc(\op,\op)$. But ${\sf K}_{\omega}$ is $\sigma$-compact and so has the Hurewicz property, and it is countable dimensional, so has property $\Sc(\op,\op)$. 

\begin{corollary}\label{extendedclass} Let $X$ be a separable metric space which has an ${\sf F}_{\sigma}$ subset $C$ such that: $C$ has the Hurewicz property, and for every open set $U\subset X$ with $C\subset U$, $X\setminus U$ is finite dimensional. Then the following are equivalent:
\begin{enumerate}
\item{$X$ has the property $\Sc(\op,\op)$.}
\item{$X$ has the property $\Sc(\op_{fin},\op)$.}
\end{enumerate}
\end{corollary}
The proof uses the fact that $\Sc(\op_{fin},\op)$ and $\Sc(\op,\op)$ are preserved by ${\sf F}_{\sigma}$-subsets.  

\section{Remarks}

In \cite{wh1} Hurewicz introduced a property weaker than the Hurewicz property, and known as Menger's property: For each sequence $(\mathcal{U}_n:n<\infty)$ of open covers of a space $X$ there is a sequence $(\mathcal{V}_n:n<\infty)$ of finite sets such that for each $n$, $\mathcal{V}_n\subset\mathcal{U}_n$, and $\bigcup_{n<\infty}\mathcal{V}_n$ is a cover of $X$.  Theorem \ref{scfinandhaver} shows that if a metrizable space has the Hurewicz property and also $\Sc(\op_{fin},\op)$, then it has the Haver property. We have the following conjecture:
\begin{conjecture}\label{mengerhaver}
There is a metrizable space $X$ with the Menger property and $\Sc(\op_{fin},\op)$, which does not have the Haver property in some metric.
\end{conjecture}
Note that Conjecture \ref{mengerhaver} implies that the answer to Borst's Question 3.10 is ``no". 

We also expect that for each $n>1$ that $\Sc(\op_n,\op)$ the implication $\Sc(\op_n,\op)\Rightarrow\Sc(\op_{n+1},\op)$ is false. 

In Remark D of \cite{erpol2} E. and R. Pol showed that a metrizable space has the property $\Sc(\op,\op)$ if, and only if, it has the Haver property in all equivalent metrics. This gives another way to conclude Theorem \ref{hurewiczscfin}: By Theorems \ref{hurewicztotbdd} and \ref{scfinandhaver}, we see that the Hurewicz property and $\Sc(\op_{fin},\op)$ implies the Haver property for all equivalent metrics. Also: By \cite{erpol} Remark D, Conjecture \ref{mengerhaver} translates to statement that Theorem \ref{hurewiczscfin} fails if the Hurewicz property is replaced with the Menger property.

\section{Acknowledgement}

I would like to thank Elzbieta and Roman Pol for communicating to me the inspiring results in \cite{erpol} and \cite{erpol2}.

\end{document}